\documentclass[11 pt,dvipsnames]{article}
\usepackage{harvard}
\usepackage{graphicx}
\usepackage{camnum}
\usepackage{amsmath}
\usepackage{amsfonts}
\usepackage{epstopdf}
\usepackage{xcolor}

\numberwithin{figure}{section}
\numberwithin{equation}{section}
\numberwithin{table}{section}
\numberwithin{theorem}{section}

\allowdisplaybreaks

\newcommand{\DDD}{\mathcal{D}}

\newcommand{\ee}{{\mathrm e}}
\newcommand{\ii}{{\mathrm i}}

\newcommand{\LLL}{\mathcal{L}}

\date{}

\title{Mathematical foundations of spectral methods for time-dependent PDEs}

\begin{document}
\thispagestyle{empty}

\author{Arieh Iserles\\
Department of Applied Mathematics and Theoretical Physics\\
Centre for Mathematical Sciences\\
University of Cambridge\\
Cambridge CB3 0WA, United Kingdom\\[4pt]
Email: \texttt{ai@maths.cam.ac.uk}}

\maketitle

\tableofcontents

\begin{abstract}
  The contention of this paper is that a spectral method for time-dependent PDEs is basically no more than a choice of an orthonormal basis of the underlying Hilbert space. This choice is governed by a long list of considerations: stability, speed of convergence, geometric numerical integration, fast approximation and efficient linear algebra. We subject different choices of orthonormal bases, focussing on the real line, to these considerations. While nothing is likely to improve upon a Fourier basis in the presence of periodic boundary conditions, the situation is considerably more interesting in other settings. We introduce two kinds of orthonormal bases, T-systems and W-systems, and investigate in detail their features. T-systems are designed to work with Cauchy boundary conditions, while W-systems are suited to zero Dirichlet boundary conditions. 
\end{abstract}

\section{Introduction: what is a spectral method?}

Spectral methods have been around since the pioneering book of \citeasnoun{gottlieb77nas} and they have been the subject of several monographs \cite{boyd01cfs,canuto06sm,funaro92pad,hesthaven07smt,trefethen00smm}. A natural setting for such methods are pure boundary-value PDEs and indeed most -- but definitely not all -- texts on spectral methods have been concerned with them. The purpose of the current paper is to revisit the theory of spectral methods for {\em time-dependent PDEs\/} and recast it in a form consistent with the emerging theory of {\em geometric numerical integration\/} \cite{blanes16cig,hairer06gni,iserles18wgn}.

Our starting point is a  time-dependent problem
\begin{equation}
  \label{eq:1.1}
  \frac{\partial u}{\partial t}=\LLL u+f,\qquad t\geq0,\quad x\in\Omega,
\end{equation}
given in a separable Hilbert space $\mathcal{H}$, equipped with the inner product $\langle\,\cdot\,,\,\cdot\,\rangle$, with suitable initial and boundary conditions. Here $\LLL$ is a well-posed differential operator. 

A spectral method is determined by an orthonormal basis $\Phi=\{\varphi_n\}_{n\in\mathbb{Z}_+}$ of $\mathcal{H}$. The solution of \R{eq:1.1} is expressed in the underlying basis,
\begin{equation}
  \label{eq:1.2}
  u(x,t)=\sum_{n=0}^\infty \hat{u}_n(t)\varphi_n(x),
\end{equation}
where the values of $\hat{u}_n(0)$ are given by expanding the initial condition in the underlying basis, while the subsequent evolution of the coefficients $\hat{u}_n(t)$ is given by the {\em Galerkin conditions\/}
\begin{equation}
  \label{eq:1.3}
  \left\langle \sum_{m=0}^\infty \hat{u}_m'(t)\varphi_m-\LLL\!\left(\sum_{m=0}^\infty \hat{u}_m(t)\varphi_m\!\right)-f,\varphi_n\!\right\rangle=0,\qquad n\in\mathbb{Z}_+.
\end{equation}
In particular, once $\LLL$ is linear, \R{eq:1.3} reduces to 
\begin{equation}
  \label{eq:1.4}
  \hat{u}_n'=\sum_{m=0}^\infty  \langle\LLL\varphi_m,\varphi_n\rangle\hat{u}_m+\langle f,\varphi_n\rangle,\qquad n\in\mathbb{Z}_+.
\end{equation}
by virtue of orthogonality. In practice, \R{eq:1.3} or \R{eq:1.4} are truncated and the subsequent ODE system is solved by a suitable method. 

And that's it! A spectral method is completely determined by an orthonormal basis of the underlying Hilbert space. In other words, the many features that determine the quality of a spectral method hinge on the choice of $\Phi$. All the rest are technicalities which, needless to say, are hugely important and nontrivial.

In the sequel we reduce the discussion for simplicity's sake to the linear case. The attributes of a spectral method that govern the choice of $\Phi$ are
\begin{description}
\item[Stability:] The truncated ODE systems
\begin{equation}
  \label{eq:1.5}
  \hat{u}_n'=\sum_{m=0}^{N}  \langle\LLL\varphi_m,\varphi_n\rangle\hat{u}_m+\langle f,\varphi_n\rangle,\qquad n=0,\ldots,N.
\end{equation}
of \R{eq:1.4} must be uniformly well posed as $N\rightarrow\infty$. This, by the agency of the Lax equivalence theorem, guarantees the convergence of finite-dimensional solutions to the solution of \R{eq:1.1} and is a non-negotiable feature of a numerical method. 

\item[Speed of convergence:] The main reason for the efficiency -- and, simultaneously, for the limitations -- of spectral methods is the following classical result. Consider the Hilbert space $\CC{L}_2((a,b),\D\mu)$ of functions defined in the real interval $(a,b)$ and square-integrable there with a smooth Borel measure $\D\mu$. Let $\Omega$ be an open domain such that $[a,b]\subset\Omega$ and suppose that both $f$ and the $\varphi_n$s are {\em analytic\/} in $\Omega$. Then the speed of convergence, once $f$ is approximated by $\Phi$, is at least exponential, in the sense that
\begin{displaymath}
  \limsup_{N\rightarrow\infty} \left\|f-\sum_{n=0}^N \hat{f}_n \varphi_n\right\|^{1/N}=\rho
\end{displaymath}
for some $\rho\in[0,1)$ \cite{trefethen13ata,walsh65iar}. This can be extended to multivariate domains. In other words, as long as both the solution of \R{eq:1.1} and the orthonormal set $\Phi$ are analytic (or, at least, very regular), spectral methods converge very rapidly. Moreover, we require fairly small $N$ for suitably accurate solution. 

In other words, under suitable conditions nothing can beat a spectral method but, once either the solution or the orthonormal basis exhibit poor regularity, we can expect very slow convergence, large $N$ and a spectral method is no longer competitive. 

\item[Geometric numerical integration:] A fairly recent paradigm in numerical analysis is {\em Geometric Numerical Integration (GNI):\/} it is not enough to produce an approximate solution to given accuracy but also it is vital to preserve qualitative features of the solution \cite{blanes16cig,hairer06gni}. In particular, two structural features lend themselves to appropriately designed spectral methods: preservation of dissipativity in the solution of the diffusion equation and preservation of mass (essentially, of the $\CC{L}_2$ norm) in the solution of the linear Schr\"odinger equation. Naive discretisation methods may fail to respect these features, thereby mis-computing fundamental physical features.

Note incidentally that preservation of dissipativity or mass automatically (and trivially) implies stability. 

\item[Fast expansion:] In the context of spectral methods, a function $f\in\mathcal{H}$ is rendered as a sequence $\hat{\MM{f}}=\{\hat{f}_n\}_{n\in\mathbb{Z}_+}$ of its expansion coefficients. In other words, we have a {\em fundamental map\/} $\mathcal{F}$ such that $\hat{\MM{f}}=\mathcal{F}f$. This map is invertible. The computation of $\mathcal{F}f$ and $\mathcal{F}^{-1}\hat{\MM{f}}$ need be done at least once per time step and it is important to do so by fast algorithms.  This is an important reason why spectral methods are so popular in the presence of periodic boundary conditions: it is then natural to choose $\Phi$ as a Fourier basis, whereby $\mathcal{F}$ and its inverse can be computed to exponential accuracy by Fast Fourier Transform.

Note that we can assume without loss of generality that $\mathcal{F}$ is unitary, $\mathcal{F}^{-1}=\mathcal{F}^*$, and an isometry -- all it takes is the right scaling of the $\ell_2$ norm. This assumption, which means that the Plancherel theorem is valid in this setting, is henceforth taken for granted.

\item[Linear algebra:] Ultimately, computer calculations (as distinct from the mathematical theory underpinning algorithms) reduce in the main to linear algebra.  While finite differences and finite elements typically yield very large, yet very sparse linear systems, spectral methods lead to considerably smaller systems (in particular for highly regular solutions, cf.\ our comments on rapid convergence requiring small $N$ to attain stipulated accuracy). In traditional spectral methods, though, these systems are dense (and often poorly conditioned). However, it turns out that in the context of time-dependent problems an appropriate choice of $\Phi$ leads to `nice' linear systems, which can be solved rapidly. Again, everything depends on the choice of an orthonormal basis!
\end{description}

To recap, these five features determine how good a spectral method is -- and the only weapon in our arsenal is the choice of an orthonormal basis. 

Before we proceed further, we need to introduce a fundamental concept.

\vspace{6pt}
\noindent{\bf Definition} Given an orthonormal basis $\Phi\subset\mathcal{H}\cap \CC{C}^1$, a {\em differentiation matrix\/} $\DDD$ is a linear map taking the vector $\MM{\varphi}=\{\varphi_n\}_{n\in\mathbb{Z}_+}$ to its derivative $\MM{\varphi}'=\{\varphi_n'\}_{n\in\mathbb{Z}_+}$. In other words,
\begin{displaymath}
  \varphi_m'=\sum_{n=0}^\infty \DDD_{m,n}\varphi_n,\qquad m\in\mathbb{Z}_+.
\end{displaymath}

A differentiation matrix always exists and is unique because $\Phi$ is a basis of $\mathcal{H}$. Using $\DDD$ in tandem with the fundamental map $\mathcal{F}$ allows us to assemble a spectral method. For example, consider (in a single space dimension) the {\em diffusion equation\/}
\begin{equation}
  \label{eq:1.6}
  \frac{\partial u}{\partial t}=\frac{\partial^2 u}{\partial x^2},\qquad x\in[-1,1],\quad t\geq0,
\end{equation}
with zero Dirichlet boundary conditions and the {\em linear Schr\"odinger equation\/} 
\begin{equation}
  \label{eq:1.7}
  \ii \frac{\partial u}{\partial t}=\frac12\frac{\partial^2 u}{\partial x^2}-V(x)u,\qquad x\in\mathbb{R},\quad t\geq0.
\end{equation}

For \R{eq:1.6} we choose an orthonormal basis $\Phi$ of \hspace*{-4pt}\raisebox{9pt}{\hspace*{6pt}\small$\circ$\hspace*{-6pt}}$\CC{H}_2^1[0,1]$ and represent the solution in the form \R{eq:1.2}, where $\hat{\MM{u}}_n(0)=\mathcal{F} u(\,\cdot\,,0)$. It follows from \R{eq:1.4} that
\begin{displaymath}
  \hat{\MM{u}}'=\DDD^2\hat{\MM{u}},\qquad \mbox{hence}\qquad \hat{\MM{u}}(t)=\ee^{t\DDD^2}\hat{\MM{u}}(0).
\end{displaymath}
In particular, once $\DDD$ is a skew-Hermitian matrix, $\DDD^2=-\DDD^\top\DDD$ is a Hermitian semi-negative matrix and we conclude that $\|\hat{\MM{u}}(t)\|_2\leq \|\hat{\MM{u}}(0)\|_2$: this implies not just stability but also dissipativity, a major structural feature of the diffusion equation.

Likewise, for \R{eq:1.7}, we obtain
\begin{displaymath}
 \hat{\MM{u}}'=\ii \left(-\frac12 \DDD^2+\mathcal{F}V\mathcal{F}^{-1}\right)\!\MM{\hat{u}}
\end{displaymath}
In particular, once $\DDD$ is skew-Hermitian and recalling the unitarity of $\mathcal{F}$, the matrix $\mathcal{E}:=-\frac12 \DDD^2+\mathcal{F}V\mathcal{F}^{-1}=\frac12 \DDD^\top\DDD+\mathcal{F}V\mathcal{F}^*$ is Hermitian and 
\begin{displaymath}
  \hat{\MM{u}}(t)=\ee^{\ii t\mathcal{E}}\hat{\MM{u}}(0)\qquad \Rightarrow\qquad \|\hat{\MM{u}}(t)\|_2=\|\hat{\MM{u}}(0)\|_2.
\end{displaymath}
Thus, we have both stability and unitarity of the solution operator: the latter is crucial to applications of \R{eq:1.7} in quantum mechanics.

A common denominator to both examples is that we wish the differentiation matrix to be skew Hermitian (or skew symmetric if it is real). This should come as little surprise because the differentiation operator (subject to Cauchy, periodic or zero Dirichlet boundary conditions) is itself skew symmetric.

The two above examples illustrate the major requirement in the remainder of the paper: {\em the differentiation matrix $\DDD$ is skew Hermitian.\/} This is a condition that will guide us in the choice of the orthonormal set $\Phi$.

We do not consider in this paper a different approach to the stability of spectral methods which, in effect, replaces the Galerkin conditions \R{eq:1.3} by {\em Galerkin--Petrov conditions\/} \cite{townsend15asp}. This provides an alternative avenue toward stable spectral methods.

\section{Orthogonal systems}

\subsection{Periodic boundary conditions}

Once \R{eq:1.1} is given in a finite interval with periodic boundary conditions, there is an obvious choice of an orthonormal system: the standard Fourier basis. In particular, in the interval $[-\pi,\pi]$ the Fourier basis is $\varphi_n(x)=\ee^{\ii nx}/\sqrt{2\pi}$, $n\in\mathbb{Z}$ -- note that $\Phi=\{\varphi_n\}_{n\in\mathbb{Z}}$ is indexed over all of $\mathbb{Z}$ but this should present no difficulties. The differentiation matrix is both diagonal and skew-Hermitian,
\begin{displaymath}
  \DDD_{m,n}=
  \begin{case}
    \ii n, & m=n,\\[4pt]
    0, & \mbox{otherwise.}
  \end{case}
\end{displaymath}
An important advantage of a Fourier basis is that the first $N$ expansion coefficients of $f\in\CC{L}_2[\pi,\pi]\cap\CC{C}^{\mathrm{per}}[\pi,\pi]$ can be approximated to very high accuracy (exponential accuracy for analytic $f$) using the {\em Fast Fourier Transform (FFT)\/} in $\O{n\log_2n}$ operations.

\subsection{Cauchy boundary conditions}

Once $\Omega=\mathbb{R}$ in \R{eq:1.1}, we impose an additional condition on $\Phi$: the differentiation matrix should be both skew-Hermitian and tridiagonal. In other words, 
\begin{equation}
  \label{eq:2.1}
  \varphi_n'=-b_{n-1}\varphi_{n-1}+\ii c_n\varphi_n+\bar{b}_n\varphi_{n+1},\qquad n\in\mathbb{Z}_+,
\end{equation}
where $b_n\in\mathbb{C}\setminus\{0\}$, $c_n\in\mathbb{R}$, $n\in\mathbb{Z}_+$, and $b_{-1}=0$. However, letting $\theta_0=0$, $\theta_{n+1}=\theta_n-\arg b_n$ and $\varphi_n=\ee^{\ii\theta_n}\tilde{\varphi}_n$, $n\in\mathbb{Z}_+$, we obtain 
\begin{displaymath}
  \tilde{\varphi}_n'=-|b_{n-1}|\tilde{\varphi}_{n-1}+\ii c_n \tilde{\varphi}_n+|b_n|\tilde{\varphi}_{n+1},\qquad n\in\mathbb{Z}_+.
\end{displaymath}
Therefore, we may assume without loss of generality that $b_n>0$, $n\in\mathbb{Z}_+$, in \R{eq:2.1}. Such an orthonormal basis $\Phi$ of $\CC{L}_2(\mathbb{R})$ is called a {\em T-system\/} and we henceforth proceed to characterise such systems. The theory of T-systems has been initiated in \cite{iserles19oss}.

Let $\vartheta_n(\xi)=(-\ii)^n\hat{\varphi}_n(\xi)$, $n\in\mathbb{Z}_+$, where $\hat{\varphi}_n$ is the Fourier transform of $\varphi_n$. Then it follows from \R{eq:2.1} that the $\vartheta_n$s obey the three-term recurrence relation
\begin{equation}
  \label{eq:2.2}
  b_n\vartheta(\xi)=(\xi-c_n)\vartheta_n(\xi)-b_{n-1}\vartheta_{n-1}(\xi),\qquad n\in\mathbb{Z}_+.
\end{equation}

It is known even to casual students of mathematical analysis that any system of orthonormal polynomials obeys \R{eq:2.2}, but a much stronger result is the {\em Favard theorem:\/} Whenever $\vartheta_0$ is a constant (and, consequently, each $\vartheta_n$ a polynomial of degree $n$), there exists a Borel measure $\D\mu$ such that $\{\vartheta_n\}_{n\in\mathbb{Z}_+}$ is an orthonormal polynomial system with respect to the inner product
\begin{displaymath}
  \langle f,g\rangle_\mu=\int_{-\infty}^\infty f(\xi)g(\xi)\D\mu(\xi)
\end{displaymath}
\cite{chihara78iop}. We conclude that the general solution of \R{eq:2.2} is $\vartheta_n=gp_n$, $n\in\mathbb{Z}_+$, where $\{p_n\}_{n\in\mathbb{Z}_+}$ is an orthonormal polynomial system with respect to some Borel measure $\D\mu$, while $g\in\CC{L}_2(\mathbb{R})$ is a nonzero distribution. Applying inverse Fourier transform to $\ii^n \vartheta_n$ results in $\Phi=\{\varphi_n\}_{n\in\mathbb{Z}_+}$ which is guaranteed to obey \R{eq:2.1} and this, indeed, characterises all its solutions. 

Yet, there is absolutely no reason {\em per se\/} why $\Phi$ should be orthonormal: for this we need to choose the function $g$ judiciously. Because of the Plancherel theorem,
\begin{displaymath}
  \int_{-\infty}^\infty \varphi_n(x)\overline{\varphi_n(x)}\D x=\int_{-\infty}^\infty \hat{\varphi}_m(\xi)\overline{\hat{\varphi}_n(\xi)}\D\xi=\int_{-\infty}^\infty p_m(\xi)q_n(\xi) |g(\xi)|^2\D\xi.
\end{displaymath}
Therefore we have orthogonality if and only if $|g(\xi)|^2\D\xi=\D\mu(\xi)$: an obvious choice is $g(\xi)=\sqrt{w(\xi)}$, where $w$ is the Radon--Nikodym derivative of $\mu$, $\D\mu(\xi)=w(\xi)\D\xi$. This is a complete characterisation of all $\CC{L}_2$-orthonormal sequences $\Phi$ that obey \R{eq:2.1}.

Finally, we need $\Phi$ to be not just orthonormal with a tridiagonal, skew-Hermitian differentiation matrix: in addition it needs also to be a basis of $\CC{L}_2(\mathbb{R})$. In general, this need not be the case. It follows from our construction that the $\CC{L}_2$ closure of $\Phi$ is the {\em Paley--Wiener space\/} $\mathcal{PW}_{\Xi}$ of all $\CC{L}_2(\mathbb{R})\cap\CC{C}(\mathbb{R})$ functions whose Fourier transform lives in $\Xi=\CC{supp}\,\D\mu$. Thus, we need $\D\mu$ to be supported in the entire real line except that, as we will see, there are ways round it.

Our first example of a T-system is  the {\em Legendre measure\/} $\chi_{(-1,1)}(\xi)\D\xi$, whereby $p_n(\xi)=\CC{P}_n(\xi)/\sqrt{n+\frac12}$, where $\CC{P}_n$ is the standard Legendre polynomial and $g(\xi)\equiv1$. The outcome is 
\begin{equation}
  \label{eq:2.3}
  \varphi_n(x)=\frac{\sqrt{n+\frac12}}{x} \CC{J}_{n+\frac12}(x),\qquad n\in\mathbb{Z}_+,
\end{equation}
spherical Bessel functions. Since the measure is supported in $(-1,1)$, they form an orthonormal basis of $\mathcal{PW}_{(-1,1)}$ and are of little interest in the construction of spectral methods. 

Next we consider the {\em Hermite measure\/} $\D\mu(\xi)=\ee^{-\xi^2}\D\xi$. The orthonormal polynomial system consists of scaled Hermite polynomials $\CC{H}_n(\xi)/\sqrt{2^nn!\sqrt{\pi}}$, $g(\xi)=\ee^{-\xi^2/2}$ and
\begin{equation}
  \label{eq:2.4}
  \varphi_n(x)=\frac{\CC{H}_n(x)}{\sqrt{2^nn!\sqrt{\pi}}} \ee^{-x^2/2},\qquad n\in\mathbb{Z}_+,
\end{equation}
Since the Hermite measure is supported on the entire real line, \R{eq:2.4} is indeed a basis of $\CC{L}_2(\mathbb{R})$. It is well known in theoretical physics, since the $\varphi_n$s are eigenfunctions of the free Schr\"odinger operator. A T-system of the special form $\varphi_n(x)=G(x)q_n(x)$, $n\in\mathbb{Z}_+$, where $G$ is a function and $q_n$ an $n$th degree polynomial, occurs exclusively with the Hermite measure \cite{iserles19oss}.

Our final example are the {\em Malmquist--Takenaka functions\/} \cite{iserles20for}. We commence from the {\em Laguerre measure\/} $\D\mu(\xi)=\ee^{-\xi}\chi_{[0,\infty)}(\xi)$, associated with Laguerre polynomials $\CC{L}_n$, hence $g(\xi)=\ee^{-\xi/2}$. The closure of $\{\varphi_n\}_{n\in\mathbb{Z}_+}$, though, is $\mathcal{PW}_{[0,\infty)}$ and,  to extend  to $\CC{L}_2(\mathbb{R})$, we combine it with its `mirror weight' $\ee^{\xi}\chi_{(-\infty,0]}(\xi)$, which corresponds to $p_n(\xi)=\CC{L}_n(-\xi)$. It is convenient to index the ensuing functions by integers,
\begin{equation}
  \label{eq:2.5}
  \varphi_n(x)=\sqrt{\frac{2}{\pi}} \ii^n \frac{(1+2\ii x)^n}{(1-2\ii x)^{n+1}},\qquad n\in\mathbb{Z}.
\end{equation}
This is a complex-valued system and so is $\DDD$, with $b_n=n+1$ and $c_n=2n+1$, $n\in\mathbb{Z}$.

A particularly appealing feature of the T-system \R{eq:2.5} is that, changing the variable $x=\tan(\theta/2)$,
\begin{displaymath}
  \hat{f}_n=\int_{-\infty}^\infty f(x)\overline{\varphi_n(x)}\D x=\frac{(-\ii)^n}{\sqrt{2\pi}} \int_{-\pi}^\pi \!\!\left(1\!-\!\ii\tan\frac{\theta}{2}\right)\! f\!\left(\frac12 \tan\frac{\theta}{2}\right)\!\ee^{-\ii n\theta}\D\theta,\quad\!\! n\in\mathbb{Z},
\end{displaymath}
a system of Fourier integrals which can be evaluated rapidly with FFT. A very minor generalisation of this setting leads to all possible T-systems that can be computed by FFT \cite{iserles20for}, while all T-systems that can be evaluated by Fast Cosine and Sine Transforms have been characterised in \cite{iserles21fco}.

\subsection{Dirichlet boundary conditions}

Can T-systems be extended to the setting of a compact interval $(a,b)$ with zero Dirichlet boundary conditions? While their construction via Fourier transform no longer works, in a forthcoming publication we introduce an alternative, more algebraic means to `build' a T-system inductively. Given $\varphi_0\in\CC{C}^\infty$, we construct $\{\varphi_n\}_{n\in\mathbb{N}}$, as well as $b_n,c_n$ by recursion for $n\in\mathbb{Z}_+$ \cite{iserles24ost}, and this works in any interval with zero Dirichlet boundary conditions. However, the problem is that once $\varphi_n(a)=0$ for all $n$ at some point $a$ then $\MM{\varphi}'=\DDD\MM{\varphi}$ implies that $\varphi_n'(a)=0$, $n\in\mathbb{Z}_+$ and, more generally, $\MM{\varphi}^{(\ell)}=\DDD^\ell\MM{\varphi}$ for all $\ell\in\mathbb{Z}_+$. This and the fact that $\DDD^\ell$ is bounded imply that $\varphi_n^{(\ell)}(a)=0$ for all $n,\ell\in\mathbb{Z}_+$. Therefore the $\varphi_n$s cannot be analytic at the endpoints. 

In principle, this can be overcome by placing essential singularities at the endpoints while maintaining integrability (and analyticity) in $(a,b)$: for example, we might take $(a,b)=(-1,1)$ and 
\begin{displaymath}
  \varphi_0(x)=a_0 \exp\!\left(-\frac{1}{1-x^2}\right)\!,\qquad a_0=\left[\int_{-1}^1 \exp\!\left(-\frac{2}{1-y^2}\right)\!\D y\right]^{-1/2},
\end{displaymath}
This can be extended to a full T-system using the methodology of \cite{iserles24ost}. Unfortunately,  essential singularities at the endpoints strike back: as $n$ grows, the $\varphi_n$s develop a boundary layer of increasing amplitude and the approximating power of $\Phi$ is very poor indeed. 

An alternative to T-systems, better suited to this setting, are W-systems \cite{iserles24ssm}.  Given $-\infty\leq a<b\leq \infty$ and a weight function $w\in\CC{C}^1(a,b)$, we let
\begin{equation}
  \label{eq:2.6}
  \varphi_n(x)=\sqrt{w(x)}p_n(x),\qquad n\in\mathbb{Z}_+,
\end{equation}
where $\{p_n\}_{n\in\mathbb{Z}_+}$ is an orthonormal polynomial system with respect to the Borel measure $\D\mu(x)=w(x)\chi_{(a,b)}(x)\D x$. Trivially, $\Phi$ is an orthonormal basis of $\CC{L}_2(a,b)$. Moreover, since
\begin{displaymath}
  \DDD_{m,n}+\DDD_{n,m}=\int_a^b (\varphi_m'\varphi_n+\varphi_m\varphi_n')\D x=\int_a^b \frac{\D}{\D x} \varphi_m \varphi_n \D x=\varphi_m(x)\varphi_n(x)\,\rule[-5pt]{0.75pt}{18pt}_{\,x=a}^{\,x=b},
\end{displaymath}
it follows from \R{eq:2.6} that $\DDD$ is skew symmetric if and only if $w(a)=w(b)=0$. 

Unlike in the case of T-systems, the differentiation matrix is neither tridiagonal nor can we take for granted that all its powers are bounded -- indeed, the entire purpose of the exercise is to construct systems so that $\DDD^{s+1}$ blows up for some $s\in\mathbb{Z}_+$: provided that $\DDD^\ell$ is bounded for $\ell\leq s$, we say that $s$ is the {\em index\/} of $w$ and denote it by $\CC{ind}\,w$. It is possible to prove that $\CC{ind}\,w\geq 2$ if and only if ${w'}^2/w$ is a weight function (i.e.\ nonnegative, with all its moments bounded and the zeroth moment positive) \cite{iserles24ssm}. Moreover, a necessary condition for $\CC{ind}\,w\geq s$ in a finite interval $(a,b)$ is that $w(x)=(x-a)^\alpha(b-n)^\beta \tilde{w}(x)$, where $\alpha,\beta>s-1$, with obvious extension to $(a,\infty)$ and $(-\infty,b)$  \cite{iserles24ssm}. 

Our first example of a W-system is when $\DDD$ is a banded matrix: in a sense, this is the least interesting case because all its powers are then bounded as well, which somewhat defeats the purpose of the exercise. This happens if and only if $w(x)=\ee^{-q(x)}$, $x\in\mathbb{R}$, where $q$ is an even-degree polynomial, $q(x)=x^{2r}+\mbox{l.o.t.}$. This corresponds to a {\em Freud measure\/} and a special case, $q(x)=x^2$, yields Hermite functions -- the only instance of a W-system being also a T-system (cf.\ \R{eq:2.4}).

Without loss of generality, it is enough to consider the intervals $[0,\infty)$ and $[-1,1]$. In $[0,\infty)$ the obvious choice is the {\em Laguerre weight function\/} $w(x)=x^\alpha \ee^{-x}$, $\alpha>-1$, whereby
\begin{equation}
  \label{eq:2.7}
  \varphi_n(x)=\sqrt{\frac{n!}{\CC{\Gamma}(n+1+\alpha)}} x^{\alpha/2} \ee^{-x/2} \CC{L}_n^{(\alpha)}(x),\qquad x>0,\quad n\in\mathbb{Z}_+.
\end{equation}
It takes several pages of dense algebra to derive the explicit form of the differentiation matrix for $\alpha>1$,
\begin{equation}
  \label{eq:2.8}
  \DDD_{m,n}=
  \begin{case}
    -\frac12 \GG{a}_m\GG{b}_n, & m\geq n+1,\\[2pt]
    0 & m=n,\\[2pt]
    \frac12 \GG{a}_n\GG{b}_m, & m\leq n-1,
  \end{case}
\end{equation}
where
\begin{displaymath}
  \GG{a}_m=\sqrt{\frac{m!}{\CC{\Gamma}(m+1+\alpha)}},\qquad \GG{b}_n=\sqrt{\frac{\CC{\Gamma}(n+1+\alpha)}{n!}}
\end{displaymath}
\cite{iserles24ssm}.  We already know that $\CC{ind}\,w<\alpha+1$ and a long combinatorial proof shows that $\DDD^\ell$ is bounded if and only if $\ell<\alpha+1$. Therefore $\CC{ind}\,w=\lfloor \alpha-1\rfloor+2$, $\alpha>1$. 

In the interval $(-1,1)$ an obvious choice is the {\em ultraspherical weight function\/} $w(x)=(1-x^2)^\alpha$, $\alpha>-1$, leading to
\begin{equation}
  \label{eq:2.9}
  \varphi_n(x)=\sqrt{\frac{n!(2n+2\alpha+1)\CC{\Gamma}(n+2\alpha+1)}{2^{2\alpha+1}[\CC{\Gamma}(n+\alpha+1)]^2}}(1-x^2)^{\alpha/2}\CC{P}_n^{(\alpha,\alpha)}(x),\qquad n\in\mathbb{Z}_+,
\end{equation}
where $\alpha>-1$. The differentiation matrix for $\alpha>1$ is
\begin{equation}
  \label{eq:2.10}
  \DDD_{m,n}=
  \begin{case}
    0, & m+n\mbox{\ even},\\[2pt]
     -\frac12 \GG{a}_m\GG{b}_n, & m+n \mbox{\ odd},\;m\geq n+1,\\[2pt]
    \frac12 \GG{a}_n\GG{b}_m, & m+n \mbox{\ odd},\; m\leq n-1
  \end{case}
\end{equation}
with
\begin{displaymath}
  \GG{a}_m=\sqrt{\frac{m!(2m+2\alpha+1)}{2\CC{\Gamma}(m+2\alpha+1)}},\qquad \GG{b}_n=\sqrt{\frac{(2n+2\alpha+1)\CC{\Gamma}(n+2\alpha+1)}{2n!}}.
\end{displaymath}
Again, $\CC{ind}\,w=\lfloor \alpha-1\rfloor+2$, $\alpha>1$ \cite{iserles24ssm}. 

As will be apparent in the next section, the forms \R{eq:2.8} and \R{eq:2.10} of differentiation matrices are a key toward efficient algebra.

\subsection{Multivariate setting}

It is straightforward to extend Fourier expansions to the torus $\mathbb{T}^d$, $d\geq2$, while we can be generalise T-systems to $\mathbb{R}^d$, $d\geq2$ and W-systems  to parallelepipeds using tensor products. Moreover, W-systems can be extended to different multivariate geometries using the theory of orthogonal polynomials in $d$ dimensions, although each generalisation is nontrivial, not least the derivation of the differentiation matrix. In particular, \cite{gao24fss} extends W-systems to $d$-dimensional balls using a generalisation of Zernike polynomials, while a generalisation to triangles and simplexes is subject to current active research.

\section{Fast expansion and linear algebra}

The following tasks are at the heart of a successful implementation of a spectral method:
\begin{itemize}
\item The computation of  $\mathcal{F}_N$ (the fundamental map, truncated to its first $N+1$ coefficients) and its inverse $\mathcal{F}^*_N$;
\item Formation of a product of the form $\DDD_N \MM{x}$, where $\DDD_N$ is an $(N+1)\times(N+1)$ principal minor of $\DDD$ and $\MM{x}\in\mathbb{C}^{N+1}$;
\item The solution of linear algebraic systems of the form $(I-\kappa \DDD_N)\MM{y}=\MM{x}$, where $\kappa\in\mathbb{C}$ -- note that a system of the form $q(\DDD_N)\MM{y}=\MM{x}$, where $q$ is a polynomial of degree $k$, can be represented as a succession of $k$ such problems. Note further that $\DDD^*=-\DDD$ means that polynomials in both $\DDD_N$ and $\DDD^*_N$ fit into this formalism.
\item The computation of $f(t\DDD_N)\MM{x}$, where $f$ is a rational function -- in particular $f(z)=\ee^z$;
\end{itemize}
These four imperatives, combined with the speed of convergence -- a subject relegated to the next section -- are fundamental in guiding us in the choice of the underlying orthonormal basis $\Phi$.

\subsection{T-systems}

Five T-systems lend themselves to fast computation of $\mathcal{F}_N$: the Malmquist--Takenaka system \R{eq:2.5} can be computed using FFT\footnote{This can be somewhat extended to `generalised Malmquist--Takenaka system' \cite{iserles20for} but there is no apparent advantage to be gained from this.} while
\begin{displaymath}
  \varphi_n(x)=\frac{(-1)^n}{\sqrt{g_n^{(\alpha,\beta)}}} (1-\tanh x)^{(\alpha+1)/2} (1+\tanh x)^{(\beta+1)/2} \CC{P}_n^{(\alpha,\beta)}(\tanh x),\qquad n\in\mathbb{Z}_+,
\end{displaymath}
where $\CC{P}^{(\alpha,\beta)}/\sqrt{g_n^{(\alpha,\beta)}}$ is the orthonormalised $n$-th degree Jacobi polynomial, can be computed by Fast Cosine Transform for $\alpha=\beta=\pm\frac12$ and by Fast Sine Transform for $\alpha=-\beta=\pm\frac12$ \cite{iserles21fco}. Expansions in Hermite functions \R{eq:2.4} can be computed rapidly using asymptotic methods for orthogonal polynomial expansions from \cite{olver29fau}. 

Linear algebra with T-systems is expedited by the fact that $\DDD_N$ is tridiagonal. Therefore, both forming $\DDD_N\MM{x}$ and solving $(I-\kappa \DDD_N)\MM{y}=\MM{x}$ requires just $\O{N}$ operations. 

There are several options to compute (or approximate) $f(t\DDD_N)\MM{x}$ and they all take advantage of the tridiagonal form of the differentiation matrix. The simplest is to approximate $f$ by a rational function, e.g.\ a Pad\'e approximation \cite{higham08fm}. Another option is to use Krylov subspace methods \cite{hochbruck97oks}. Probably the most intriguing is to use the {\em Dunford formula\/} (known also as Cauchy--Dunford, Dunford--Taylor, Riesz and Fantappi\`e formula, \cite{caratelli21fda},
\begin{equation}
  \label{eq:3.1}
  f(T)=\frac{1}{2\pi\ii} \oint_{\gamma} f(z)(zI-T)^{-1}\D z,
\end{equation}
where $\gamma$ is a simple, closed smooth curve with positive orientation, enclosing the spectrum of the operator $T$, such that $f$ is holomorphic inside and on $\gamma$. The integral \R{eq:3.1} can be discretised as a Riemann sum at roots of unity of a parametric representation of $\gamma$ and this lends itself to fast evaluation of $f(T_N)\MM{x}$ (where $T_N$ is an $(N+1)\times(N+1)$ section of $T$) using the Fast Fourier Transform \cite{hackbusch15hma}.

\subsection{W-systems}

Insofar as W-systems are concerned, computing $\mathcal{F}_N$ is the same as expanding $f/\sqrt{w}$ in orthogonal polynomials and we can leverage the asymptotic methods from \cite{olver29fau} for  Laguerre and ultraspherical W-systems \R{eq:2.7} and \R{eq:2.9}, respectively.

A differentiation matrix of a W-system is dense (unless the weight function is a Freud weight) and we no longer enjoy the benefit of tridiagonality. However, linear algebra with both Laguerre and ultraspherical W-systems can be accomplished very effectively exploiting another feature of their differentiation matrices: they are {\em semiseparable\/} \cite{hackbusch15hma,vandebril08mc1,vandebril08mc2}. Specifically, a matrix is semiseparable of rank $r$ if each finite-dimensional sub-matrix which lies entirely above or entirely underneath the main diagonal  is at most of rank $r\in\mathbb{N}$ and there exists at least one such matrix of rank exactly  $r$. It follows at once from \R{eq:2.8} and  \R{eq:2.10} that $r=1$ for both Laguerre and ultraspherical W-systems.

As an example how semiseparability leads to fast linear algebra, we show how to compute rapidly the product $\DDD_N\MM{x}$ for the Laguerre W-system, while noting that the solution of linear systems (and much more) can be also done in an economical manner using the material in \cite{vandebril08mc1}. Actually, we can do better: compute 
\begin{displaymath}
  y_m=\sum_{n=0}^K \DDD_{m,n}x_n,\qquad m=0,\ldots,N,
\end{displaymath}
where $N\leq K$, as this allows for far greater accuracy in truncating $\DDD$. To this end, recalling \R{eq:2.8}, we define
\begin{displaymath}
  \sigma_m=\sum_{n=0}^{m-1} \GG{b}_n x_n,\quad \rho_m=\sum_{n=m+1}^K \GG{a}_n f_n,\qquad m=0,\ldots,N.
\end{displaymath}
Then
\begin{displaymath}
  y_0=\GG{b}_0\rho_0,\qquad y_m=-\GG{a}_m\sigma_m+\GG{b}_m\rho_m,\quad m=1,\ldots,N.
\end{displaymath}
Since
\begin{displaymath}
  \sigma_m=\sigma_{m-1}+\GG{b}_{m-1}x_{m-1},\qquad \rho_m=\rho_{m-1}-\GG{a}_m x_m,
\end{displaymath}
once $\{\GG{a}_n\}_{n=0}^K$ and $\{\GG{b}_n\}_{n=0}^K$ are precomputed, it takes just $K+5N$ flops to evaluate $\MM{y}$, while naive computation `costs' $\O{NK}$ flops.

\section{Speed of convergence}

It is convenient to commence from W-systems in a compact interval, specifically considering the ultraspherical W-system in $(-1,1)$ (although other W-systems in a compact interval are subject to the same ground rules). The speed of convergence is governed by standard theory of orthogonal polynomials. Thus, let $\rho>0$ and consider the {\em Bernstein ellipse\/}
\begin{displaymath}
  \mathcal{B}_\rho=\{\rho\ee^{\ii\theta}+\rho^{-1}\ee^{-\ii\theta}\,:\, -\pi\leq\theta\leq\pi\}.
\end{displaymath}
Suppose that $f$ is analytic inside $\mathcal{B}_\rho$ and let $\{p_n\}_{n\in\mathbb{Z}_+}$ be an orthonormal polynomial system in $(-1,1)$ with respect to the measure $\D\mu$. Then 
\begin{equation}
  \label{eq:4.1}
  \limsup_{n\rightarrow\infty} \left\|f-\sum_{m=0}^n \check{f}_n p_n\right\|_{\CC{L}_2[(-1,1),\D\mu]}^{1/n}=\frac{1}{\rho},
\end{equation}
where 
\begin{displaymath}
  \check{f}_n=\int_{-1}^1 f(x) p_n(x)\D\mu(x),\qquad n\in\mathbb{Z}_+.
\end{displaymath}
In other words, the speed of convergence is exponential and governed by the  `analyticity ellipse' \cite{trefethen00smm,walsh65iar}.

Let $f$ be analytic inside $\mathcal{B}_\rho$ and, in addition, $f^{(i)}(-1)=f^{(i)}(1)=0$, $i=0,\ldots,s-1$. Therefore,
\begin{displaymath}
  f^{\langle s\rangle}(x):=\frac{f(x)}{(1-x^2)^{s}}
\end{displaymath}
is analytic inside $\mathcal{B}_\rho$. Recalling that $\varphi_n(x)=(1-x^2)^{\alpha/2} \tilde{\CC{P}}_n^{(\alpha,\alpha)}(x)$, where $\tilde{\CC{P}}_n^{(\alpha,\alpha)}$ is an orthonormal Jacobi polynomial, we note that it is analytic if and only if $\alpha$ is an even nonnegative integer. Moreover, we have 
\begin{displaymath}
  \hat{f}_n=\int_{-1}^1f(x) \varphi_n(x)\D x=\int_{-1}^1 (1-x^2)^{\alpha/2+s} f^{\langle s\rangle}(x) \tilde{\CC{P}}_n^{(\alpha,\alpha)}(x)\D x,\qquad n\in\mathbb{Z}_+.
\end{displaymath}
Therefore, once we choose $\alpha=2s$, it is true that $\hat{f}_n=\check{f}_n^{\langle s\rangle}$ and \R{eq:4.1} implies exponential convergence. Note in passing that any $\alpha=2q$, where $q\geq s$ is an integer, guarantees exponential convergence in the $\CC{L}_2$ norm but, once $q>s$, $\CC{H}_\infty^s$ convergence is not possible because the $(s+1)$st derivative of the approximation -- but not of $f$ -- vanishes at the endpoints. Thus, the choice $\alpha=2s$ is the `sweet spot', allowing for both exponential convergence in $\CC{L}_2$ and for $\CC{H}_\infty^s$ convergence.

\begin{figure}[htb]
  \begin{center}
    \includegraphics[width=180pt]{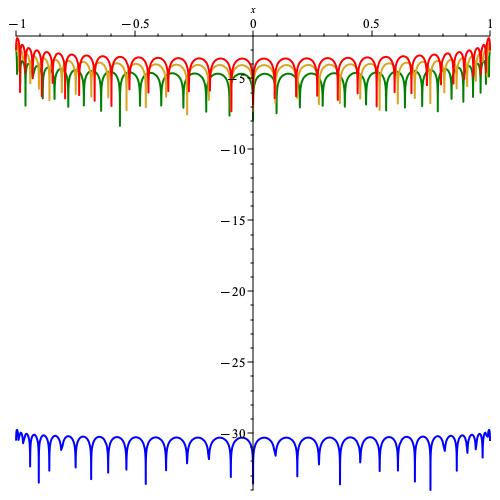}\hspace*{15pt}\includegraphics[width=180pt]{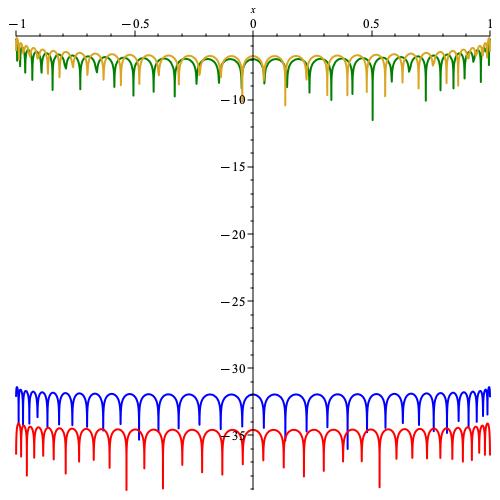}
    \caption{Pointwise errors in $[-1,1]$ for ultraspherical W-systems and the functions $f(x)=\sin\pi x$ (on the left) and $f(x)=\cos^2(\pi x/2)$. In each case $=\alpha=1,2,3,4$ corresponds to green, blue, yellow and red lines.}
    \label{fig:4.1}
  \end{center}
\end{figure}

We provide two examples for ultraspherical W-systems. Firstly, we take $f(x)=\sin\pi x$ and display pointwise (i.e.\ $\CC{L}_\infty$) errors in Fig.~\ref{fig:4.1} (left) for $\alpha=1,2,3,4$ using 31 first terms in the expansion. The difference is striking: while $\alpha=1,3,4$ yield very slow convergence, the `sweet spot' $\alpha=2$ results in a  truly spectacular performance. On the right of Fig.~\ref{fig:4.1} we display the error for the function $f(x)=\cos^2(\pi x/2)$: note that both $f$ and $f'$ vanish at $\pm1$. This means that the `sweet spot' is $\alpha=4$ and this is confirmed in the plot -- although $\alpha=2$ also demonstrates good performance. 

\begin{figure}[htb]
  \begin{center}
    \includegraphics[width=180pt]{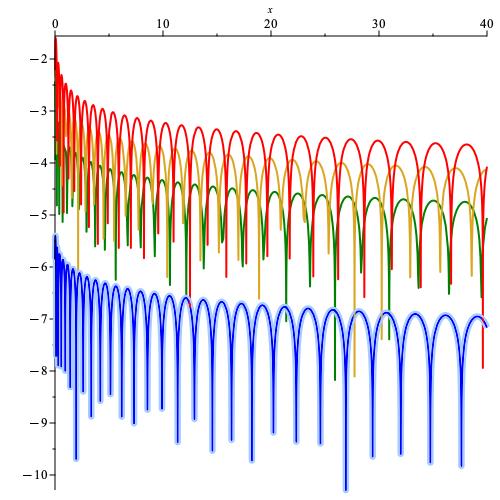}\hspace*{15pt}\includegraphics[width=180pt]{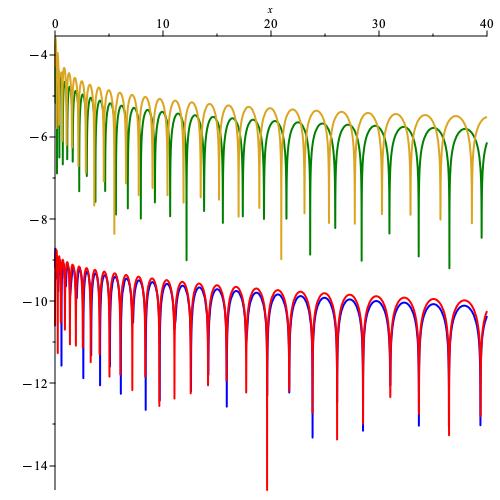}
    \caption{Pointwise errors in $[-0,\infty)$ for Laguerre W-systems and the functions $f(x)=x\ee^{-x}/(1+x)$ (on the left) and $f(x)=x\ee^{-2x}\sin x$. In each case $=\alpha=1,2,3,4$ corresponds to green, blue, yellow and red lines.}
    \label{fig:4.2}
  \end{center}
\end{figure}

In Fig.~\ref{fig:4.2} we turn our attention to Laguerre W-systems, using an identical argument in choosing optimal $\alpha$: on the left  the function $f(x)=x\ee^{-x}/(1+x)$ (note that $f(0)=0$, $f'(0)\neq0$) and on the right $f(x)=x\ee^{-2x}\sin x$: in that case both $f$ and $f'$ vanish at the origin. Again, the `sweet spot' $\alpha=2$ performs significantly better on the left, while $\alpha=2$ and $\alpha=4$, both corresponding to analytic $\varphi_n$s, demonstrate roughly similar behaviour, far superior to $\alpha=1$ and $\alpha=3$. The difference is much less striking, though, than in Fig.~\ref{fig:4.1} -- for a good reason. Unlike in compact integrals, there is no comprehensive convergence theory of analytic orthonormal bases approximating analytic functions in $[0,\infty)$. The standard Bernstein ellipse argument is no longer valid and there is no alternative theory. While there is a body of evidence to suggest fractional exponential convergence, e.g.\ at the rate of $c\ee^{-d n^\gamma}$ for $\gamma\in(0,1)$, this falls short of complete theory. 

The situation is, if at all, even more open-ended with regard to T-systems, because comprehensive convergence theory on the real line is not available either. There are few results available and they are intriguing, sometimes baffling. In one case, when $f\in\CC{L}_2(\mathbb{R})$ is a rational function, everything is known and can be derived with relative ease. Thus, suppose that the poles of $f$ are $s_1,\ldots,s_N$ and that $\Im s_k<0$ for $k=1,\ldots,K-1$ and $\Im s_k>0$ for $k=K,\ldots,N$ -- clearly, unless $\Im s_k\neq0$, $k=1,\ldots,N$, the function $f$ cannot be $\CC{L}_2(\mathbb{R})$! Set
\begin{displaymath}
  \sigma_k=\frac{1-2\ii s_k}{2s_k-\ii},\qquad k=1,\ldots,N
\end{displaymath}
and note that $|\sigma_k|<1$ for $k=1,\ldots,K-1$ and $|\sigma_k|>1$ for $k=K,\ldots,N$. Then
\begin{equation}
  \label{eq:4.2}
  \limsup_{n\rightarrow\infty}\left\|f-\sum_{k=-n}^n \hat{f}_n \varphi_n\right\|^{1/n}_{\CC{L}_2(\mathbb{R})}=\max\left\{\max_{k=1,\ldots,K-1}|\sigma_k|,\max_{k=K,\ldots,N} |\sigma_k|^{-1}\right\}<1.
\end{equation}
We have exponential convergence! Unfortunately, all this goes away even upon a seemingly minor change of scenery. Thus, consider first $f(x)=1/(1+x+x^2)$: we have $s_1=-\frac12-\ii\frac{\sqrt{3}}{2}$, $s_2=-\frac12+\ii\frac{\sqrt{3}}{2}$, therefore $\sigma_1=-(2+3\ii)/(5+2\sqrt{3})$ and $\sigma_2=-(2+3\ii)/(5-2\sqrt{3})$ -- note that $|\sigma_1|=|\sigma_2|^{-1}=\sqrt{(37-20\sqrt{3})/13}\approx 0.4260$ and all is well. However, suppose that we replace $f$ by $\cos x/(1+x+x^2)$: on the face of it, not much has changed -- the cosine is an entire function, bounded by 1 in modulus along the real line. However, it has an essential singularity at infinity (i.e.\ at the North Pole of the Riemann sphere) and this degrades the speed of convergence down to $\O{|n|^{-4/3}}$ -- from an exponential heaven to little better than linear convergence.

The situation is, if at all, worse for a Hermite T-system: now both $1/(1+x+x^2)$ and $\cos x/(1+x+x^2)$ exhibit $\O{n^{4/3}}$ convergence. 

The speed of convergence for different flavours of T-systems is an open problem. Asymptotic analysis in \cite{weideman95cht} provides few clues, as does the analysis of Hermite and Malmquist--Takenaka systems in the case of wave packets in \cite{iserles23awp}, but much is yet to be done to investigate how well T-systems approximate analytic (or  rougher) functions.

\bibliographystyle{agsm}

\end{document}